\documentclass[12pt]{amsart}

\usepackage{amsmath}
\usepackage{amssymb}
\usepackage{amscd}
\usepackage{indentfirst}
\usepackage[T1]{fontenc}
\usepackage[frenchb]{babel}

\theoremstyle{plain}
\newtheorem{thm}{Théorème}[section]
\newtheorem{lem}[thm]{Lemme}
\newtheorem{cor}[thm]{Corollaire}
\newtheorem{prop}[thm]{Proposition}
\theoremstyle{definition}

\newtheorem{ntt}[thm]{Notations}
\theoremstyle{remark}
\newtheorem{rmq}[thm]{Remarque}
\newtheorem{exm}[thm]{Exemple}

\DeclareMathOperator{\pic}{Pic}
\DeclareMathOperator{\disc}{disc}
\newcommand{\A}{\mathcal{A}}
\newcommand{\B}{\mathcal{B}}
\newcommand{\Ao}{\mathcal{A}^{\circ}}

\newcommand{\E}{\mathcal{E}}
\newcommand{\Eo}{\mathcal{E}^{\circ}}

\newcommand{\gm}{\mathbf{G}_{{\rm m}}}
\newcommand{\gmk}{\mathbf{G}_{{\rm m},K}}
\newcommand{\gmt}{\mathbf{G}_{{\rm m},T}}

\newcommand{\limp}{\varprojlim\,}

\newcommand{\gabis}{^\varphi\Gamma'}
\newcommand{\gan}{^n\Gamma'}

\newcommand{\dai}{\mathcal{D}}
\newcommand{\e}{h_*}
\newcommand{\He}{\mathcal{H}}

\DeclareMathOperator{\homr}{Hom}

\newcommand{\Hom}{\underline{\rm Hom}}
\newcommand{\Ext}{\underline{\rm Ext}^1}

\newcommand{\ext}{{\rm Ext}^1}

\begin{document}

\title{Variétés abéliennes et invariants arithmétiques}

\author{Jean Gillibert}

\address{
Université de Caen,
Laboratoire de Mathématiques Nicolas Oresme (CNRS umr 6139),
BP 5186, 14032 Caen cedex, France}

\email{jean.gillibert@math.unicaen.fr}

\date{31 Août 2005}

\maketitle


\section{Introduction}

Soient $R$ un anneau de Dedekind excellent, de corps de fractions $K$, et
$S=Spec(R)$.
On
considère
un $S$-schéma en groupes commutatif
$G$, fini et plat sur $S$, et l'on note $G^D$ le dual de Cartier de $G$.
Nous disposons
d'un
homomorphisme
\[
\begin{CD}
\pi:H^1(S,G^D) @>>> \pic(G)\\
\end{CD}
\]
explicité en premier par Waterhouse (voir \cite{w}, Theorem 5). Si l'on
considère que la
notion
de
torseur (sous un schéma en groupes
fini) généralise celle d'extension galoisienne, alors on peut dire que
$\pi$ mesure la
structure
galoisienne des $G^D$-torseurs, le
groupe $\pic(G)$ étant identifié à un groupe de classes (voir \cite{cnt}).

Plusieurs auteurs se sont intéressés à la construction de $G^D$-torseurs
dont l'image
par $\pi$
est
triviale, c'est-à-dire de torseurs
dont la structure galoisienne est triviale. Plus précisément, supposons
que $G$ soit un
sous-groupe
d'un $S$-schéma abélien $A$. Soit
$B:=A/G$ le schéma abélien quotient, et soit $A^t$ (resp. $B^t$) le schéma
abélien dual
de $A$
(resp. $B$). Alors (par dualité) nous avons une suite exacte (de faisceaux
abéliens pour la topologie fppf sur
$S$)
\[
\begin{CD}
0 @>>> G^D @>>> B^t @>>> A^t @>>> 0\\
\end{CD}
\]
qui donne lieu, par application du foncteur des sections globales, à un
morphisme cobord
$\partial:A^t(S)\rightarrow H^1(S,G^D)$. Nous
obtenons ainsi un homomorphisme $\psi:=\pi\circ\partial$, introduit en
premier (dans le
cas où
$G=A[n]$ et où $R$ est l'anneau des
entiers d'un corps de nombres) par M. J. Taylor \cite{t2}, et que l'on
appelle
usuellement {\em
class-invariant homomorphism}. Srivastav
et Taylor \cite{st}, puis Agboola \cite{a2}, et enfin Pappas \cite{p1} ont
montré que,
si $A$ est
une courbe elliptique et si l'ordre de
$G$ est premier à $6$, alors les points de torsion de $A^t(S)$
appartiennent au noyau de
$\psi$.

Dans \cite{gil}, l'auteur a généralisé la construction de $\psi$ dans le
cas où $G$ est
un
sous-groupe fini et plat d'un $S$-schéma en
groupes semi-stable, ainsi que le résultat d'annulation sur les points de
torsion (dans
le cas
d'une courbe elliptique semi-stable, en
supposant l'ordre de $G$ premier à $6$).

Notre but est d'étendre cette construction dans le cadre suivant : soient
$\A$ le modèle
de Néron
d'une $K$-variété abélienne
semi-stable, $\Ao$ la composante neutre de $\A$, $\Gamma$ un sous-groupe
({\em i.e.} un
sous-schéma
en groupes ouvert) de $\Phi:=\A/\Ao$
(quotient pour la topologie fppf sur $S$) et $G$ un sous-groupe fermé,
quasi-fini et
plat de
$\A^{\Gamma}$ (où $\A^{\Gamma}$ désigne
l'image réciproque de $\Gamma$ dans $\A$). Nous définissons (paragraphe
\ref{defpsiz})
un
homomorphisme $\psi$ (associé à l'inclusion
$G\subseteq\A^{\Gamma}$) qui se factorise de la façon suivante
\[
\begin{CD}
\A^{t,\Gamma'}(S) @>>> H^1(S,\Hom_S(G,\gm)) @>\pi>> \pic(G)\,,\\
\end{CD}
\]
où $\Gamma'$ est l'orthogonal de $\Gamma$ sous l'accouplement de
monodromie défini par
Grothendieck
(voir le paragraphe \ref{neron}), et
où $\pi$ est un homomorphisme qui généralise celui de Waterhouse (voir le
paragraphe
\ref{picard}).

Signalons ici que l'application $\Gamma\mapsto\A^{\Gamma}$ réalise une
bijection entre
l'ensemble
des sous-groupes ouverts de $\Phi$ et
l'ensemble des modèles semi-stables de $\A_K$.

Du point de vue technique, la principale différence avec la situation
considérée dans
\cite{gil}
est l'absence de dual de Cartier
naturel pour le schéma en groupes $G$. Cependant, on donne ici une preuve
de la nullité
du
faisceau
$\Ext_S(G,\gm)$ pour la topologie
fppf (voir le lemme \ref{youpi}). Ce résultat, qui nous a été communiqué
par L.
Moret-Bailly,
nous
permet de généraliser de façon
naturelle les constructions précédentes.

Travailler avec des groupes quasi-finis permet de reformuler le résultat
d'annulation de
$\psi$
démontré dans \cite{gil}. Ainsi on peut
énoncer, comme corollaire du théorème 4.1 de \cite{gil}, le résultat
suivant (voir le
paragraphe
\ref{annulsurtors}) :

\begin{thm}
\label{intro1}
Soit $\E$ le modèle de Néron d'une courbe elliptique à réduction
semi-stable sur $K$, et
soit
$m>0$
un entier naturel premier à $6$.
Alors les homomorphismes
\[
\begin{CD}
\E^t(S) @>>> H^1(S,\Hom_S(\Eo[m],\gm)) @>\pi>> \pic(\Eo[m])\\
\end{CD}
\]
(associé à $\Eo[m]\subseteq \Eo$), et
\[
\begin{CD}
\E^{t,\circ}(S) @>>> H^1(S,\Hom_S(\E[m],\gm)) @>\pi>> \pic(\E[m])\\
\end{CD}
\]
(associé à $\E[m]\subseteq \E$) s'annulent sur les points de torsion.
\end{thm}

\begin{rmq}
On constate que les groupes $\Eo[m]$ et $\E[m]$ sont affines (paragraphe
\ref{affine}).
D'autre
part, si les points de $m$-torsion de
$\E_{\eta}$ sont tous $K$-rationnels, alors le groupe $\E[m]$ est fini et
plat sur $S$
(\cite{gil},
prop. 3.6), donc
$\Hom_S(\E[m],\gm)=\E[m]^D$ est le dual de Cartier de $\E[m]$. Par contre,
si $\E$ n'a
pas
partout
bonne réduction, $\Eo[m]$ n'est
jamais un groupe fini, sauf pour $m=1$, valeur pour laquelle il est nul
(voir la
remarque
\ref{boxall}).
\end{rmq}

L'un des objectifs de ce travail est de << passer à la limite >> dans
l'étude de $\psi$.
Supposons
que $K$ soit un corps de nombres, et
que $R$ soit l'anneau des entiers de $K$. Sous ces hypothèses, nous
introduisons (cf.
section 4)
une version arakélovienne $\hat\psi$ de
notre homomorphisme (généralisant celle introduite par Agboola et Pappas
dans \cite{ap})
: soit
$\ell$ un nombre premier, on note alors
$$\hat{\Psi}_{\ell}=\limp\hat{\psi}_{\ell^n}:\A^t(S)\otimes\mathbf{Z}_{\ell}\longrightarrow
\limp
\widehat{\pic}(\Ao[\ell^n])$$
la flèche obtenue par passage à la limite projective de l'homomorphisme
$\hat{\psi}_{\ell^n}$
associé à l'inclusion
$\Ao[\ell^n]\subseteq\Ao$ (cf. paragraphe \ref{inj}). Soient
$\disc(K/\mathbb{Q})$ le
discriminant de $K/\mathbb{Q}$, et $U\subseteq S$ l'ouvert de bonne
réduction de $\A$.
Les
résultats de \cite{ap} impliquent alors
le résultat qui suit :

\begin{thm}
\label{intro2}
L'homomorphisme $\hat{\Psi}_{\ell}$ est injectif modulo les points de
torsion. En outre,
si tous
les points de $S$ de caractéristique $\ell$ sont contenus dans $U$, et si
$\ell$ ne
divise pas
$6\cdot\disc(K/\mathbb{Q})$, alors $\hat{\Psi}_{\ell}$ est injectif.
\end{thm}

\begin{rmq}
Supposons que, pour toute place $v$ de mauvaise réduction de $\A$, $\ell$
ne divise pas
l'ordre
de
$\Phi_v$. Alors nous avons
$\Ao[\ell^n]=\A[\ell^n]$ pour tout entier $n$, ce qui permet d'exprimer
plus simplement
l'ensemble
d'arrivée de $\hat{\Psi}_{\ell}$.
\end{rmq}

\begin{rmq}
Philippe Cassou-Noguès et Martin Taylor ont remarqué une analogie entre
l'annulation de $\psi$ sur les points de torsion et la conjecture de Birch
et Swinnerton-Dyer (voir la remarque 5.8 de \cite{t3} ainsi que les
commentaires qui suivent le théorème 4 de \cite{cnt}).
Soit à présent $\Psi_{\ell}:\A^t(S)\otimes\mathbf{Z}_{\ell}\rightarrow
\limp \pic(\Ao[\ell^n])$ le morphisme déduit de $\hat{\Psi}_{\ell}$ par
oubli des métriques. Les remarques précédentes suggèrent un lien entre
l'injectivité de $\Psi_{\ell}$ sur les points d'ordre infini et la
conjecture de Birch et Swinerton-Dyer $\ell$-adique.

Remarquons que, si $\A$ est une courbe elliptique à multiplication
complexe, ayant partout bonne réduction sur $K$ et réduction ordinaire en
$\ell$, alors l'injectivité de $\Psi_{\ell}$ modulo les points de torsion
a été démontrée (sous certaines hypothèses) par Agboola et Taylor (voir
\cite{at} ou le théorème 6 de \cite{cnt}).
\end{rmq}

\noindent\textit{Remerciements.} Je remercie Laurent Moret-Bailly pour m'avoir communiqué la preuve du lemme \ref{youpi}, à la lumière duquel cet article a été largement remanié, ainsi que pour sa relecture de l'ensemble du texte. Je tiens également à exprimer ma reconnaissance envers John Boxall pour l'encadrement de ma thèse, dont ce travail est issu.


\section{Problèmes d'exactitude}

Rappelons les notations qui seront en vigueur tout au long de cet article
: $R$ est un
anneau de
Dedekind excellent, de corps de
fractions $K$. Soit $S=Spec(R)$ et soit $\eta$ le point générique de $S$.
Nous noterons
$\gm$ le
groupe multiplicatif sur $S$.

On dit qu'un $S$-schéma en groupes est semi-stable s'il est commutatif,
lisse, séparé,
et si les
composantes neutres de ses fibres sont
extensions de variétés abéliennes par des tores.

On fixe une fois pour toutes un $S$-schéma en groupes semi-stable $A$ dont
la fibre
générique
$A_{\eta}$ est une variété abélienne. On
notera $U\subseteq S$ l'ouvert de bonne réduction de $A$, de sorte que
$A_U$ est un
$U$-schéma
abélien.

\subsection{Isogénies duales sur le petit site fppf}
\label{isoduz}
Supposons que l'on se donne un épimorphisme (fppf) $f:A\rightarrow B$
entre $S$-schémas
en
groupes
semi-stables, tel que
$f_{\eta}:A_{\eta}\rightarrow B_{\eta}$ soit une isogénie. Alors $\ker f$
est un
$S$-schéma en
groupes plat et quasi-fini (voir
\cite{ray}, \S\/7.3, lemma 1). De plus, nous avons un carré cartésien
\[
\begin{CD}
\ker f @>>> A \\
@VVV @VVfV \\
S @>e>> B \\
\end{CD}
\]
où $e:S\rightarrow B$ désigne la section unité, qui est une immersion
fermée car $B$ est
séparé.
Par suite, $\ker f\rightarrow A$ est
une immersion fermée.

Réciproquement, soit $G$ un sous-$S$-schéma en groupes fermé, quasi-fini
et plat de $A$.
Le lemme
suivant montre que $G$ s'inscrit dans
une suite exacte.

\begin{lem}
Le faisceau quotient $A/G$, pour la topologie fppf sur $S$, est
représentable par un
$S$-schéma
en
groupes semi-stable.
\end{lem}

\begin{proof}
Le faisceau quotient $A/G$ est représentable (\cite{anan}, chap. IV,
théorème 4.C) par
un
$S$-schéma en groupes, que nous noterons $B$
(on se sert du fait que $A$ est de type fini sur $S$, et $S$ régulier de
dimension $\leq
1$). La
projection canonique
$\varphi:A\rightarrow B$ est fidèlement plate, et $A$ est un $S$-schéma
plat, donc $B$
est un
$S$-schéma plat. De plus $B$ est lisse sur
$S$ grâce au critère de lissité par fibres (\cite{ray}, \S\/2.4, prop. 8).
Enfin, les
composantes
neutres des fibres de $B$ sont
extensions de variétés abéliennes par des tores.
\end{proof}

Ainsi nous avons une suite exacte
\[
\begin{CD}
0 @>>> G @>>> A @>\varphi>> B @>>> 0\\
\end{CD}
\]
de faisceaux abéliens pour la topologie fppf sur $S$. De plus, la
restriction
$\varphi_U:A_U\rightarrow B_U$ est une isogénie entre
$U$-schémas abéliens. Il existe alors une isogénie duale
$\varphi_U^t:B_U^t\rightarrow
A_U^t$,
dont
le noyau est le dual de Cartier
$G_U^D$ de $G_U$. Ici, $A_U^t$ et $B_U^t$ désignent les $U$-schémas
abéliens duaux de
$A_U$ et
$B_U$ respectivement.

Nous voulons prolonger l'isogénie duale $\varphi_U^t$ en un morphisme de
faisceaux sur
$S$, pour
cela nous allons nous servir du
foncteur $\Hom(-,\gm)$ et de ses dérivés. Dans cette optique, le gros site
fppf ne nous
convient
pas, car le faisceau $\Hom(A,\gm)$
n'est pas nul. Nous allons donc utiliser un autre site, qui nous permettra
d'énoncer le
lemme
\ref{nullz}.

Plus précisément, nous considérons le << petit site fppf >> sur $S$ (resp.
sur $U$),
c'est-à-dire
la catégorie des schémas plats sur $S$
(resp. sur $U$) munie d'une structure de site pour la topologie fppf. Nous
appellerons
faisceau
(pour la topologie fppf) sur $S$ (resp.
sur $U$) un faisceau sur ce site.

Nous noterons $\Hom_S$ (resp. $\Hom_U$) le faisceau des homomorphismes de
faisceaux
abéliens,
restreint au petit site fppf sur $S$
(resp. sur $U$).

Soit $j:U\rightarrow S$ l'inclusion, et soit $j^*$ le foncteur << image
inverse >> de
faisceaux
correspondant. La flèche $j:U\rightarrow
S$ étant un objet du petit site fppf sur $S$, le foncteur $j^*$ est un <<
foncteur de
localisation
>> (voir \cite{gro4}, exposé IV,
paragraphes 5.1 à 5.4).

En particulier, il en résulte que l'image par $j^*$ d'un faisceau
représentable (disons
par un
$S$-schéma plat $Y$) est représenté par
le $U$-schéma $Y_U:=Y\times_S U$.
D'autre part, si $F_1$ et $F_2$ sont deux faisceaux abéliens sur $S$, la
flèche
canonique
$$j^*(\Hom_S(F_1,F_2))\rightarrow \Hom_U(j^*F_1,j^*F_2)$$
est un isomorphisme (voir \cite{gro4}, exposé IV, prop. 12.3, b), p. 502).
De plus,
$j^*$ est
exact
et admet un adjoint à gauche $j_{!}$
exact. Par suite, $j^*$ envoie les injectifs sur des injectifs (voir
\cite{gro4}, exposé
V,
paragraphe 2.2). Nous dérivons alors des
deux côtés de la flèche et obtenons un isomorphisme
$j^*(\Ext_S(F_1,F_2))\simeq
\Ext_U(j^*F_1,j^*F_2)$. En particulier, soit $Y$ un
$S$-schéma en groupes plat tel que $Y_U$ soit un $U$-schéma abélien, alors
nous obtenons
un
isomorphisme
\begin{equation}
\label{iminvz}
j^*(\Ext_S(Y,\gm))\simeq Y_U^t
\end{equation}
où $Y_U^t$ est le schéma abélien dual de $Y_U$. On se sert ici du fait que
$\Ext_U(Y_U,\gm)$ est
isomorphe à $Y_U^t$ (voir \cite{gro7},
exposé VII, 2.9.5 et 2.9.6).

Nous avons le résultat suivant (voir \cite{gil}, lemme 2.2) :

\begin{lem}
\label{nullz}
Soit $Y$ un $S$-schéma en groupes plat, dont la fibre générique $Y_{\eta}$
est une
variété
abélienne. Alors $\Hom_S(Y,\gm)$ est nul.
\end{lem}

D'autre part, nous devons à L. Moret-Bailly la démonstration du lemme qui
suit (laquelle
dépend
du
fait que $S$ est régulier de
dimension $\leq 1$).

\begin{lem}[L. Moret-Bailly]
\label{youpi}
Soit $G$ un schéma en groupes (commutatif) plat, séparé et quasi-fini sur
$S$. Alors
$\Ext_S(G,\gm)$ est nul.
\end{lem}

\begin{proof}
On peut supposer $S$ local hensélien. Le groupe $G$ étant quasi-fini et
séparé sur $S$,
il admet
un
plus grand sous-groupe ouvert et
fermé $H$ fini sur $S$ (voir \cite{gro7}, exposé IX, 2.2.3). De plus, la
platitude de
$G$
entraîne
celle de $H$. Le quotient $G/H$ est
alors étale sur $S$ (sa section unité est ouverte), de fibre spéciale
nulle. Sa fibre
générique
est
donc un $K$-schéma en groupes
(automatiquement séparé) fini étale. On en déduit que $G/H$ est lui-même
plat,
quasi-fini et
séparé
sur $S$. Le groupe $H$ étant fini et
plat sur $S$, nous avons $\Ext_S(H,\gm)=0$ (voir \cite{gro7}, exposé VIII,
3.3.1). Ainsi
on se
ramène au cas où $G$ est étale, à fibre
spéciale nulle.

Soit $\Omega$ une extension de $G$ par $\gm$. On remarque que $G_K$ et la
section unité
de $G$
forment un recouvrement ouvert de $G$.
Donc trivialiser l'extension $\Omega$ équivaut à trivialiser sa fibre
générique
(extension de
$G_K$
par $\gmk$). Une telle
trivialisation existe après extension finie de $K$, donc après revêtement
fini et plat
de $S$,
d'où
le résultat.
\end{proof}

A présent, nous considérons la suite
\[
\begin{CD}
0 @>>> G @>>> A @>\varphi>> B @>>> 0\\
\end{CD}
\]
comme étant une suite exacte de faisceaux sur le petit site fppf de $S$.

Nous obtenons, en appliquant le foncteur $\Hom_S(-,\gm)$ à cette suite,
une (longue)
suite exacte
de cohomologie
\begin{align*}
\Hom_S(A,\gm) \rightarrow \Hom_S(G,\gm) \rightarrow \Ext_S(B,\gm)
\rightarrow
\Ext_S(A,\gm)
\rightarrow \Ext_S(G,\gm)
\end{align*}
dont les termes sont des faisceaux abéliens sur $S$. Le premier terme est
nul, d'après
le lemme
\ref{nullz}, ainsi que le dernier terme,
d'après le lemme \ref{youpi}. Ainsi nous obtenons une suite exacte
\begin{equation}
\label{suite}
\begin{CD}
0 \longrightarrow \Hom_S(G,\gm) \longrightarrow \Ext_S(B,\gm)
@>\varphi^*>>
\Ext_S(A,\gm)
\longrightarrow 0\,.\\
\end{CD}
\end{equation}
D'après ce qui précède (voir (\ref{iminvz})), son image par le foncteur
$j^*$ est la
suite
\[
\begin{CD}
0 @>>> G_U^D @>>> B_U^t @>\varphi^t_U>> A_U^t @>>> 0\,.\\
\end{CD}
\]
Cette dernière admet donc un << prolongement >> sur $S$.

Par application du foncteur des sections globales, nous pouvons déduire de
la suite
(\ref{suite})
un morphisme cobord
$$\delta:\Ext_S(A,\gm)\rightarrow H^1(S,\Hom_S(G,\gm))\,.$$

Nous allons voir à présent comment la théorie des biextensions permet
d'éclaircir les
choses en
construisant explicitement des sections
du faisceau $\Ext_S(A,\gm)$.

\subsection{Biextensions et faisceau image}
\label{neron}
Soit $\A$ le modèle de Néron de $A_{\eta}$. Alors, $A$ étant semi-stable,
la flèche
$A\rightarrow
\A$ prolongeant l'application
identique $A_{\eta}\rightarrow A_{\eta}$ est une immersion ouverte, et
induit un
isomorphisme
entre
les composantes neutres (voir
\cite{ray}, \S\/7.4, prop. 3). Par suite, nous identifierons $A$ à un
sous-groupe ouvert
de $\A$.
Soit $\Phi:=\A/\Ao$ le groupe des
composantes de $\A$. Alors il existe un sous-faisceau $\Gamma$ de $\Phi$
tel que $A$
soit l'image
réciproque de $\Gamma$ par la
surjection canonique $\A\rightarrow \Phi$. Nous adopterons donc la
notation usuelle
$A=\A^{\Gamma}$.

Soit à présent $A_{\eta}^t$ la variété abélienne duale de $A_{\eta}$, et
soit $\A^t$ son
modèle
de
Néron, alors $(\A^t)_U=A_U^t$ est le
schéma abélien dual de $A_U$, ce qui est consistant avec les notations
précédentes. Nous
allons
voir comment la théorie des biextensions
permet d'établir un lien entre le faisceau $\Ext_S(A,\gm)$ et le schéma
$\A^t$.

On sait que la dualité entre $A_{\eta}$ et $A_{\eta}^t$ découle de
l'existence d'un
fibré en
droites $\mathcal{P}_{\eta}$ sur
$A_{\eta}\times_K A_{\eta}^t$, que l'on appelle fibré de Poincaré.
Cependant nous
envisageons ici
la dualité dans un cadre plus général
à l'aide de la notion de biextension, introduite par Mumford dans
\cite{mu}. Pour une
définition
précise de cette notion nous renvoyons
le lecteur aux exposés de Grothendieck (\cite{gro7}, exposé VII) et de
Milne
(\cite{mil},
Appendix
C).

Grâce au théorème du carré, on peut munir le fibré de Poincaré
$\mathcal{P}_{\eta}$
d'une unique
structure de biextension de
$(A_{\eta},A_{\eta}^t)$ par $\gmk$, que l'on appelle la biextension de
Weil, et que l'on
note
$W_{\eta}$. Une question naturelle se pose
: $W_{\eta}$ se prolonge-t-elle en une biextension sur les modèles de
Néron ?

L'accouplement (dit << de monodromie >>) introduit par Grothendieck dans
\cite{gro7}
nous donne
la
réponse.

Plus précisément, soit $\Phi'$ le groupe des composantes de $\A^t$. On
déduit de la
lecture de
(\cite{gro7}, exposé VIII, théorème 7.1,
b)) un accouplement (associé à $W_{\eta}$)
$$\Phi\times_S\Phi'\longrightarrow ({\bf Q}/{\bf Z})_S\,.$$
Signalons ici que, $\A$ étant semi-stable, cet accouplement est non
dégénéré (voir
\cite{gro7},
exposé IX, théorème 11.5 ; on pourra
également consulter \cite{werner} pour la propriété de compatibilité
laissée au lecteur
par
Grothendieck). Le problème de prolongement
est alors résumé par la proposition qui suit.

\begin{prop}
\label{woz}
Soit $M$ (resp. $M'$) un sous-groupe de $\Phi$ (resp. $\Phi'$). Alors il
existe une
(unique)
biextension $W$ de $(\A^{M},\A^{t,M'})$ par
$\gm$ prolongeant la biextension de Weil $W_{\eta}$ sur
$(A_{\eta},A_{\eta}^t)$ si et
seulement
si
$M$ et $M'$ sont orthogonaux sous
l'accouplement.
\end{prop}

\begin{proof}
Le résultat découle de (\cite{gro7}, exposé VIII, théorème 7.1, b)) dans
le cas
particulier où la
base est un trait. En outre, après
lecture de (\cite{gro7}, exposé VIII, remarque 7.2), on voit que ce
résultat s'étend au
spectre
d'un anneau de Dedekind, ce qui est bien
le cas ici.
\end{proof}

\begin{ntt}
Soit $X$ une biextension de $(\A^{M},\A^{t,M'})$ par $\gm$. Nous noterons
$t(X)$ le
$\gm$-torseur
sous-jacent à la biextension $X$.
\end{ntt}

Nous noterons $\Gamma'$ l'orthogonal de $\Gamma$ sous l'accouplement de
monodromie, et
$W$
l'unique
biextension de
$(\A^{\Gamma},\A^{t,\Gamma'})$ par $\gm$ prolongeant $W_{\eta}$. Alors $W$
définit un
morphisme
de
faisceaux
$$\alpha:\A^{t,\Gamma'}\longrightarrow \Ext_S(\A^{\Gamma},\gm)\,.$$
Nous noterons $\gamma$ le morphisme induit par $\alpha$ sur les
$S$-sections.

On note $\B$ le modèle de Néron de la variété abélienne $B_{\eta}$, et
$\Psi$ le groupe
des
composantes de $\B$. Alors (voir \cite{ray},
\S\/7.4, prop. 3) $B$ s'identifie à un sous-groupe ouvert de $\B$. Nous
noterons donc
$B=\B^{\Lambda}$, où $\Lambda$ est le groupe des
composantes de $B$, identifié à un sous-faisceau de $\Psi$.
Soit $B_{\eta}^t$ la variété abélienne duale de $B_{\eta}$, soit $\B^t$
son modèle de
Néron, et
soit $\varphi^t:\B^t\rightarrow\A^t$
l'unique morphisme prolongeant l'isogénie duale de $\varphi_{\eta}$.

Soit $W_{\eta}'$ la biextension de Weil sur $(B_{\eta},B_{\eta}^t)$. Alors
les
biextensions
$(\varphi_{\eta}\times {\rm
id}_{B_{\eta}^t})^*(W_{\eta}')$ et $({\rm id}_{A_{\eta}}\times
\varphi_{\eta}^t)^*(W_{\eta})$
sont
isomorphes sur
$(A_{\eta},B_{\eta}^t)$.

Notons à présent $\overline\varphi:\Phi\rightarrow\Psi$ (resp.
$\overline\varphi^t:\Psi'\rightarrow\Phi'$) le morphisme induit par
$\varphi$ (resp. $\varphi^t$) sur les groupes de composantes respectifs.
On constate que
$\Lambda$
n'est autre que $\Gamma/\overline G$,
où $\overline G$ est l'image de $G$ dans $\Phi$. Par suite, nous avons
$\Lambda=\overline\varphi(\Gamma)$.
Soit le diagramme :
\[
\begin{CD}
\Phi @.\times_S @. \Phi' @>>> ({\bf Q}/{\bf Z})_S \\
@V\overline\varphi VV @. @AA\overline\varphi^t A @| \\
\Psi @.\times_S @. \Psi' @>>> ({\bf Q}/{\bf Z})_S \\
\end{CD}
\]
dans lequel on note $<,>_{\A}$ l'accouplement du haut (associé à
$W_{\eta}$), et
$<,>_{\B}$
l'accouplement du bas (associé à
$W_{\eta}'$). Il résulte alors de (\cite{gro7}, exposé VIII, 7.3.1) et de
l'identité
$(\varphi_{\eta}\times {\rm
id}_{B_{\eta}^t})^*(W_{\eta}')=({\rm id}_{A_{\eta}}\times
\varphi_{\eta}^t)^*(W_{\eta})$
que le
diagramme précédent est commutatif,
c'est-à-dire que nous avons, pour tout $x\in \Phi$ et tout $y\in \Psi'$,
l'égalité
$$<\overline\varphi(x),y>_{\B}=<x,\overline\varphi^t(y)>_{\A}\,.$$
On rappelle que $\Gamma'$ est l'orthogonal de $\Gamma$ sous l'accouplement
$<,>_{\A}$.
L'identité
précédente permet alors de montrer que
l'orthogonal de $\overline\varphi(\Gamma)$ sous l'accouplement $<,>_{\B}$
est égal à
$(\overline\varphi^t)^{-1}(\Gamma')$. Ainsi,
l'orthogonal $\Lambda'$ de $\Lambda$ est donné par :
$\Lambda'=(\overline\varphi^t)^{-1}(\Gamma')$.
Il est alors commode d'introduire
les notations suivantes :

\begin{ntt}
\label{ancdef}
On note $\gabis$ le sous-groupe de $\Phi'$ défini par
$$
\gabis=
\overline\varphi^t((\overline\varphi^t)^{-1}(\Gamma'))=\overline\varphi^t(\Lambda')
$$
de sorte que $\A^{t,\gabis}$ est l'image de $\B^{t,\Lambda'}$ par
$\varphi^t$.
\end{ntt}

\begin{rmq}
Il est clair que $\gabis$ est un sous-groupe de $\Gamma'$, donc $\Gamma$
et $\gabis$
sont
orthogonaux. Dans le cas particulier où
$\ker(\varphi^t)$ est un sous-groupe de $\A^{t,\circ}$, on a l'égalité
$\gabis=\Gamma'$.
Dans le
cas où $\Gamma'$ est nul, il est clair
que $\gabis$ l'est également. Pour un exemple où $\gabis\neq \Gamma'$,
nous renvoyons à
la
remarque
\ref{gabispapareil}.
\end{rmq}

Soit $W'$ la biextension sur $(\B^{\Lambda},\B^{t,\Lambda'})$ prolongeant
$W_{\eta}'$
dont
l'existence est assurée par la proposition
\ref{woz}, et soit $\beta:\B^{t,\Lambda'}\rightarrow
\Ext_S(\B^{\Lambda},\gm)$ le
morphisme de
faisceaux correspondant.
Pour résumer la situation, nous avons un diagramme commutatif à lignes
exactes
\begin{equation}
\label{doublesuite}
\begin{CD}
0 @>>> \ker\varphi^t @>>> \B^{t,\Lambda'} @>\varphi^t>> \A^{t,\gabis} @>>>
0 \\
@. @VVV @V\beta VV @V\alpha VV \\
0 @>>> \Hom_S(G,\gm) @>>> \Ext_S(\B^{\Lambda},\gm) @>\varphi^*>>
\Ext_S(\A^{\Gamma},\gm)
@>>> 0
\\
\end{CD}
\end{equation}
dans lequel la suite du bas n'est autre que la suite (\ref{suite}), et où
$\ker\varphi^t$ est le
noyau de
$\varphi^t:\B^t\rightarrow\A^t$, en plus d'être le noyau de la restriction
de
$\varphi^t$ à
$\B^{t,\Lambda'}$. Ce dernier fait découle
aisément de l'égalité $\Lambda'=(\overline\varphi^t)^{-1}(\Gamma')$
combinée au lemme du
serpent.

\subsection{Sur les groupes quasi-finis}
\label{affine}
En fait, notre groupe $G$ est affine. De façon plus générale, nous pouvons
énoncer la
proposition
suivante, annoncée en premier par
Raynaud (voir \cite{anan}, chap. II, prop. 2.3.1). La démonstration dépend
crucialement
du fait
que
la base $S$ est un schéma n\oe thérien
de dimension $\leq 1$.

\begin{prop}
Soit $G$ un $S$-schéma en groupes, de type fini, plat, séparé sur $S$, à
fibre générique
affine.
Alors $G$ est affine.
\end{prop}

Sous les hypothèses précédentes, nous pouvons donc écrire $G=Spec(\He)$,
où $\He$ est
une
$R$-algèbre de Hopf fidèlement plate.  Il est
clair que $\He$ est de type fini en tant que $R$-algèbre. Par contre, il
est faux en
général que
$\He$ soit de type fini en tant que
$R$-module. En effet, cette propriété équivaut au fait que $G$ soit fini
sur $S$.



\section{Invariant de Picard et homomorphisme de classes}

Nous commençons par généraliser une construction due à W. Waterhouse (voir
\cite{w},
section 2).
Puis nous définissons un homomorphisme
de classes $\psi$ qui généralise les constructions précédentes (\cite{t2},
\cite{p1},
\cite{gil}).
Nous suivons ici une démarche
semblable à celle de \cite{gil}, en abrégeant certains détails.

\subsection{L'invariant de Picard}
\label{picard}
Soit $G$ un $S$-schéma en groupes commutatif. Nous disposons d'une suite
exacte de
groupes
abéliens
(déduite de la suite spectrale
locale-globale pour les Ext (voir \cite{gro4}, exposé V, proposition 6.3,
3))
\begin{equation}
\label{exsqz}
\begin{CD}
H^1(S,\Hom_S(G,\gm)) @>\nu>> \ext(G,\gm) @>>> \Ext_S(G,\gm)(S) \\
\end{CD}
\end{equation}
où le morphisme $\nu$ est injectif. Si l'on suppose en outre que $G$ est
quasi-fini,
plat et
séparé
sur $S$, alors $\Ext_S(G,\gm)=0$
d'après le lemme \ref{youpi}, donc $\nu$ est bijectif. Nous donnons alors
la
construction
explicite
d'une flèche
$$\rho:\ext(G,\gm)\rightarrow H^1(S,\Hom_S(G,\gm))$$
telle que $\rho\circ\nu={\rm id}$. Nous obtiendrons ainsi (voir le
théorème \ref{wabis})
un
analogue, dans le cas quasi-fini, du Theorem
2' de \cite{w}.

On suppose à présent que $G$ est quasi-fini, plat et séparé sur $S$.
Définissons $\rho$
: soit
une
extension $\Omega\in\ext(G,\gm)$.
Considérons la suite exacte
\[
\begin{CD}
0 @>>> \gm @>>> \Omega @>g>> G @>>> 0\,.\\
\end{CD}
\]
Elle donne lieu, par application du foncteur $\Hom_S(G,-)$, à une suite
exacte
\[
\begin{CD}
0 @>>> \Hom_S(G,\gm) @>>> \Hom_S(G,\Omega) @>g\circ->> \Hom_S(G,G) @>>>
0\,.\\
\end{CD}
\]
On obtient ainsi un morphisme
$$\underline\delta:\homr(G,G)\rightarrow H^1(S,\Hom_S(G,\gm))\,.$$

On note alors $\rho(\Omega)$ le $\Hom_S(G,\gm)$-torseur
$\underline\delta({\rm id}_G)$.
Autrement
dit, $\rho(\Omega)$ est le faisceau
des sections $s:G\rightarrow \Omega$, au sens de la théorie des
extensions. On vérifie
qu'on
définit ainsi un morphisme de groupes
$\rho$. En résumé, nous avons :

\begin{thm}
\label{wabis}
On suppose que $G$ est quasi-fini, plat et séparé sur $S$. Alors
l'application $\rho$
définie
ci-dessus induit un isomorphisme 
$$\ext(G,\gm)\simeq H^1(S,\Hom_S(G,\gm))\,.$$
L'isomorphisme inverse $\rho^{-1}$ est égal à la flèche $\nu$ de la suite
$(\ref{exsqz})$.
\end{thm}

En composant $\nu$ avec le morphisme naturel $l:\ext(G,\gm)\rightarrow
\pic(G)$,
on obtient un homomorphisme
$$\pi:H^1(S,\Hom_S(G,\gm))\longrightarrow \pic(G)\,.$$
Dans le cas où $G$ est fini et plat, notre $\pi$ coïncide avec
l'homomorphisme défini
par
Waterhouse (voir \cite{w}, Theorem 5). En
effet, dans ce cadre, $\Hom_S(G,\gm)$ est le dual de Cartier $G^D$ de $G$.
D'autre part,
$G^D$
étant fini et plat, un argument de
descente montre que le groupe $H^1(S,G^D)$ reste inchangé, qu'il soit
calculé dans le
petit site
fppf, le gros site fppf, ou le gros
site fpqc.


\subsection{Définition et propriétés de l'homomorphisme}
\label{defpsiz}

Le résultat qui suit est démontré dans \cite{gil} comme application du
lemme
\ref{nullz}.

\begin{lem}
\label{extz}
Soit $Y$ un $S$-schéma en groupes plat, dont la fibre générique $Y_{\eta}$
est une
variété
abélienne. Alors $\Ext_S(Y,\gm)$ est
isomorphe au faisceau $T\mapsto\ext(Y_T,\gmt)$.
\end{lem}

Reprenons les notations du paragraphe \ref{neron}. Soit $E$ une extension
de
$\A^{\Gamma}$ par
$\gm$. On associe à $E$ le
$\Hom_S(G,\gm)$-torseur $\delta(E)$. Le lemme \ref{extz} permet de décrire
$\delta(E)$
comme
étant
le faisceau des extensions $\Theta$
de $\B^{\Lambda}$ par $\gm$ telles que $\varphi^*\Theta=E$.

D'autre part, en considérant le morphisme $i:G\rightarrow \A^{\Gamma}$, on
peut associer
à $E$
une
extension $i^*E$ de $G$ par $\gm$,
puis un $\Hom_S(G,\gm)$-torseur $\rho(i^*E)$. On peut décrire $\rho(i^*E)$
comme étant
le
faisceau
des sections de $i^*E$.

Or la donnée d'une section de $i^*E$ équivaut à la donnée d'une extension
$\Theta$ de
$\A^{\Gamma}$
par $\gm$ telle que
$\varphi^*\Theta=E$. On en déduit le résultat suivant :

\begin{lem}
\label{diagrefz}
Le diagramme suivant :
\[
\begin{CD}
\pic(\A^{\Gamma}) @<l^1<< \ext(\A^{\Gamma},\gm) @>\delta>>
H^1(S,\Hom_S(G,\gm))\\
@VVV @Vi^*VV @| \\
\pic(G) @<l<< \ext(G,\gm) @>\rho>> H^1(S,\Hom_S(G,\gm)) \\
\end{CD}
\]
est commutatif, où $l^1$ est le morphisme naturel
$\ext(\A^{\Gamma},\gm)\rightarrow\pic(\A^{\Gamma})$.
\end{lem}

On déduit du lemme \ref{diagrefz} le diagramme commutatif :
\begin{equation}
\label{star}
\begin{CD}
\A^{t,\Gamma'}(S) @>\gamma>> \ext(\A^{\Gamma},\gm) @>\delta>>
H^1(S,\Hom_S(G,\gm)) \\
@. @Vl^1 VV @VV\pi V \\
@. \pic(\A^{\Gamma}) @>>> \pic(G) \\
\end{CD}
\end{equation}
et on définit l'homomorphisme $\psi:\A^{t,\Gamma'}(S)\rightarrow\pic(G)$
(associé à
l'inclusion
$G\subseteq\A^{\Gamma}$) comme étant le
composé de ces morphismes.

\begin{rmq}
\label{ddef}
Soit $\dai$ l'application obtenue en composant les flèches suivantes
$$
\begin{CD}
\dai:\A^{t,\Gamma'}(S) @>\gamma>> \ext(\A^{\Gamma},\gm) @>l^1>>
\pic(\A^{\Gamma})\,.\\
\end{CD}
$$
Le diagramme (\ref{star}) montre que, pour tout $p$, $\psi(p)$ est la
restriction de
$\dai(p)$ à
$G$, que nous noterons $\dai(p)|_G$.
Ceci généralise la description géométrique de $\psi$ (obtenue par Agboola
\cite{a1} dans
le cas
où
$\A$ est un schéma abélien).

D'autre part, soit $p\in \A^{t,\Gamma'}(S)$, nous avons alors
$$\dai(p)=l^1(({\rm id}_{\A^{\Gamma}}\times p)^*(W))=({\rm
id}_{\A^{\Gamma}}\times
p)^*(t(W))\,,$$
d'où, en notant $i:G\rightarrow \A^{\Gamma}$ l'inclusion,
$$\dai(p)|_G=(i\times p)^*(t(W))\,.$$
\end{rmq}

\begin{rmq}
Supposons que $G$ soit un $S$-schéma en groupes fini et plat. Alors
l'homomorphisme
$\psi$
associé
à l'inclusion $G\subseteq \A$
coïncide avec l'homomorphisme de classes de \cite{gil}, en considérant
(avec les
notations de
\cite{gil}) les schémas en groupes
semi-stables $A=\A$ et $A'=\A^{t,\circ}$. Si de plus $G$ est contenu dans
$\Ao$, alors
l'homomorphisme $\psi$ associé à l'inclusion
$G\subseteq \Ao$ est l'homomorphisme de \cite{gil} pour $A=\Ao$ et
$A'=\A^t$.
\end{rmq}

\begin{prop}
\label{oodeux}
Soit $p\in \A^{t,\Gamma'}(S)$ et soit $N$ un entier premier à l'ordre de
$G_{\eta}$.
Alors
$\dai(p)|_G=0$ si et seulement si
$\dai(Np)|_G=0$. En particulier, $\dai(p)|_G=0$ si $p$ est un point de
$N$-torsion.
\end{prop}

\begin{proof}
Soit $m$ l'ordre de $G_{\eta}$, alors $G_{\eta}$ est tué par $m$. De même
le groupe $G$,
qui est
l'adhérence schématique de $G_{\eta}$
dans $\A^{\Gamma}$, est tué par $m$. Par suite, la multiplication par $m$
dans le groupe
$\ext(G,\gm)$, induite par la multiplication
par $m$ dans $G$, est l'application nulle. Les entiers $m$ et $N$ étant
premiers entre
eux, la
multiplication par $N$ est donc un
automorphisme du groupe $\ext(G,\gm)$. Or nous pouvons écrire
$$\dai(p)|_G=l((i\times p)^*(W))\,.$$
Autrement dit, l'homomorphisme $p\mapsto \dai(p)|_G$ se factorise à
travers le groupe
$\ext(G,\gm)$. On en déduit aisément le résultat.
\end{proof}

\subsection{Cas particulier : suite de Kummer}
\label{kumm}
Un premier avantage de notre construction est de pouvoir traiter le cas
kummérien, qui
était
l'approche originale dans \cite{t2}.

On fixe à présent un entier naturel $n>0$, et on note $n\Gamma$ l'image de
la
multiplication par
$n$ dans le groupe $\Gamma$. Alors
(voir \cite{gro7}, exposé IX, 2.2.1), $\A$ étant semi-stable, le morphisme
$[n]:\Ao\rightarrow\Ao$
est plat, surjectif, et quasi-fini.
En appliquant le lemme du serpent nous en déduisons une suite exacte
\[
\begin{CD}
0 @>>> \A^{\Gamma}[n] @>>> \A^{\Gamma} @>[n]>> \A^{n\Gamma} @>>> 0\\
\end{CD}
\]
de faisceaux abéliens pour la topologie fppf sur $S$. Il s'ensuit, d'après
les
considérations
faites au début du paragraphe
\ref{isoduz}, que $\A^{\Gamma}[n]$ est un sous-$S$-schéma en groupes
fermé, quasi-fini
et plat de
$\A^{\Gamma}$.

D'autre part, l'accouplement de monodromie étant non dégénéré,
l'orthogonal de $n\Gamma$
n'est
autre que $n^{-1}\Gamma'$, c'est-à-dire
l'image réciproque par la multiplication par $n$ du sous-groupe $\Gamma'$
de $\Phi'$. On
peut
résumer la situation par le diagramme (à
lignes exactes) :
\[
\begin{CD}
0 @>>> \A^t[n] @>>> \A^{t,n^{-1}\Gamma'} @>[n]>> \A^{t,\gan} @>>> 0 \\
@. @Vh VV @VVV  @VVV \\
0 @>>> \Hom_S(\A^{\Gamma}[n],\gm) @>>> \Ext_S(\A^{n\Gamma},\gm) @>[n]>>
\Ext_S(\A^{\Gamma},\gm)
@>>> 0 \\
\end{CD}
\]
(voir le paragraphe \ref{neron}). Ici nous avons noté $\gan$ le groupe
$n(n^{-1}\Gamma')$, en
accord avec les notations \ref{ancdef}.
Les lignes horizontales de ce diagramme sont deux prolongements de la
suite exacte de
Kummer pour
le $U$-schéma abélien
$A_U^t$.

\begin{rmq}
\label{gabispapareil}
Supposons que le groupe $\Phi'$ soit non nul et tué par l'entier $n$.
Posons $\Gamma=0$,
alors
$\Gamma'=\Phi'$ donc $\gan={}
^n\Phi'=n\Phi'=0$. Par suite, $\gan$ est distinct de $\Gamma'$. Dans
l'exemple
\ref{subtil} on
explicite également un cas pour lequel
$\A^{t,n\Phi'}(S)\neq \A^t(S)$.
\end{rmq}

Par application du foncteur des sections globales, on obtient un diagramme
commutatif
entre les
(longues) suites exactes de cohomologie
associées, qui s'écrit
\[
\begin{CD}
@>>> \A^{t,n^{-1}\Gamma'} @>[n]>> \A^{t,\gan}(S) @>\partial>>
H^1(S,\A^t[n]) @>>> \\
@. @VVV @VV \gamma V  @VV\e V \\
@>>> \ext(\A^{n\Gamma},\gm) @>[n]>> \ext(\A^{\Gamma},\gm) @>\delta>>
H^1(S,\Hom_S(\A^{\Gamma}[n],\gm))  @>>> \\
\end{CD}
\]
où $\partial$ est l'homomorphisme cobord.

La commutativité du carré à droite se traduit alors de la façon suivante :
soit $p$ dans
$\A^{t,\gan}(S)$, notons
$[n]^{-1}(p):=\partial(p)$ le $\A^t[n]$-torseur obtenu en divisant le
point $p$ par
$[n]$ dans le
faisceau $\A^t$. Alors le torseur
$(\delta\circ\gamma)(p)$ est le torseur $[n]^{-1}(p)$ auquel on a fait
subir un
changement de
groupe structural via la flèche
$h:\A^t[n]\rightarrow \Hom_S(\A^{\Gamma}[n],\gm)$.

\begin{rmq}
Si les groupes $\A^{\Gamma}[n]$ et $\A^t[n]$ sont tous les deux finis sur
$S$, alors ils
sont
duaux
l'un de l'autre au sens de Cartier,
donc $h:\A^t[n]\rightarrow \Hom_S(\A^{\Gamma}[n],\gm)$ est un
isomorphisme. On peut
alors
identifier les torseurs $[n]^{-1}(p)$ et
$(\delta\circ\gamma)(p)$. Ceci se produit en particulier lorsque $\A$ a
partout bonne
réduction.
\end{rmq}

Soit $\psi_n$ l'homomorphisme de classes associé à l'inclusion
$\A^{\Gamma}[n]\subseteq
\A^{\Gamma}$. La discussion précédente montre
que la restriction de $\psi_n$ au sous-groupe $\A^{t,\gan}(S)$ de
$\A^{t,\Gamma'}(S)$
coïncide
avec
le morphisme obtenu par composition
des flèches suivantes :
\[
\begin{CD}
\A^{t,\gan}(S) @>\partial>> H^1(S,\A^t[n]) @>\e>>
H^1(S,\Hom_S(\A^{\Gamma}[n],\gm))
@>\pi>>
\pic(\A^{\Gamma}[n])\,.\\
\end{CD}
\]
Autrement dit, pour tout $p\in \A^{t,\gan}(S)$, l'invariant $\psi_n(p)$
étudie la
structure
galoisienne du $\A^t[n]$-torseur
$[n]^{-1}(p)$ dans le groupe $\pic(\A^{\Gamma}[n])$. La seule différence
avec l'approche
originale
est que l'on fait d'abord subir au
torseur $[n]^{-1}(p)$ un changement de groupe structural.

\begin{rmq}
Dans le cas général, la suite exacte (\ref{doublesuite}) permet de donner
une
interprétation
analogue de $\psi$. Nous avons préféré nous
limiter ici au cas kummérien afin d'expliciter un peu mieux les objets
entrant en jeu.
\end{rmq}

\begin{rmq}
\label{boxall}
Supposons que $\E$ soit le modèle de Néron d'une courbe elliptique
semi-stable, n'ayant
pas
partout
bonne réduction, et soit $m>1$ un
entier naturel. Si $v$ est une place de bonne réduction, alors $\Eo_v$ est
une courbe
elliptique,
donc $\Eo_v[m]$ est un $k_v$-schéma en
groupes fini de rang $m^2$. Par contre, si $v$ est une place de mauvaise
réduction,
alors $\Eo_v$
est une forme tordue de $\gm$ sur
$k_v$, donc $\Eo_v[m]$ est un $k_v$-schéma en groupes fini de rang $m$.
Par suite,
$\Eo[m]$ n'est
pas de rang constant sur $S$, donc
n'est pas fini sur $S$.

De plus, si $m$ est premier aux ordres des groupes des composantes des
fibres de $\E$,
alors $\Eo[m]=\E[m]$, donc $\E[m]$ n'est pas fini non plus. Sous cet
éclairage, les
entiers $m$
tels que $\E[m]$ soit fini sont
exceptionnels. Signalons cependant que deux critères de finitude pour
$\E[m]$ sont
rappelés dans
le
paragraphe 3.4 de \cite{gil}.
\end{rmq}


\section{Raffinement par la théorie d'Arakelov}

Dans toute cette section, nous supposons que $K$ est un corps de nombres
et que $R$ est
l'anneau
des entiers de $K$. Nous fixons en
outre un ensemble fini $\mathfrak{S}$ de places de $K$, contenant toutes
les places
infinies de
$K$.

Soit $X$ un $S$-schéma. Un fibré inversible métrisé sur $X$ (relativement
à
$\mathfrak{S}$) est
un
fibré en droites $\mathcal{L}$ muni
d'une métrique en chacune des places de $\mathfrak{S}$. L'ensemble des
classes
d'isomorphie de
tels
objets forme un groupe (l'opération
étant induite par le produit tensoriel), noté $\widehat{\pic}(X)$. On
dispose bien sûr
d'un
homomorphisme naturel
$$\widehat{\pic}(X)\rightarrow \pic(X)$$
qui consiste à oublier les métriques sur $\mathcal{L}$.

Soit $H$ un $S$-schéma en groupes, tel que, pour toute place
$v\in\mathfrak{S}$, le
groupe
$H(\overline K_v)$ soit réunion filtrante de
sous-groupes compacts. Fixons une extension $\Omega\in\ext(H,\gm)$. Nous
avons une suite
exacte
(pour la topologie fppf) :
\[
\begin{CD}
0 @>>> \gm @>>> \Omega @>>> H @>>> 0\,.\\
\end{CD}
\]
Soit $v\in \mathfrak{S}$. Le corps $\overline K_v$ étant algébriquement
clos, il en
résulte une
suite exacte de groupes abéliens
ordinaires
\[
\begin{CD}
0 @>>> \overline K_v^* @>>> \Omega(\overline K_v) @>>> H(\overline K_v)
@>>> 0\,.\\
\end{CD}
\]
La valeur absolue $v$-adique peut-être vue comme un homomorphisme continu
$\overline
K_v^*\rightarrow \mathbb{R}$ (ici $\mathbb{R}$ est
considéré en tant que groupe additif). On déduit donc de notre extension
une extension
(de
groupes
topologiques) de $H(\overline K_v)$
par $\mathbb{R}$. Une métrique sur (le fibré inversible associé à)
$\Omega_{\overline
K_v}$ est
la
même chose qu'une trivialisation du
torseur sous-jacent à cette extension. Mais on sait que $\homr(H(\overline
K_v),
\mathbb{R})=\ext(H(\overline K_v),\mathbb{R})=0$, le
groupe $H(\overline K_v)$ étant réunion filtrante de sous-groupes
compacts. On en déduit
donc une
trivialisation canonique de
l'extension en question, d'où une métrique canonique sur
$\Omega_{\overline K_v}$.

Ainsi le fibré naturellement associé à $\Omega$ se retrouve muni d'une
métrique. On
définit ainsi
un morphisme
$$\ext(H,\gm)\rightarrow \widehat{\pic}(H)$$
qui relève l'homomorphisme naturel (nous le noterons $\hat{l}$ dans le cas
où $H=G$, et
$\hat{l}^1$
dans le cas où $H=\A^{\Gamma}$). On
obtient, en composant $\nu$ avec $\hat{l}$, un morphisme
$$\hat{\pi}:H^1(S,\Hom_S(G,\gm))\rightarrow \widehat{\pic}(G)$$
qui relève l'homomorphisme $\pi$ (défini dans le paragraphe \ref{picard}).
On pourra
consulter
(\cite{ap}, section 2) pour plus de
détails.

A l'aide du lemme \ref{diagrefz}, on vérifie que le diagramme suivant
commute :
\begin{equation*}
\begin{CD}
\A^{t,\Gamma'}(S) @>\gamma>> \ext(\A^{\Gamma},\gm) @>\delta>>
H^1(S,\Hom_S(G,\gm)) \\
@. @V\hat{l}^1 VV @VV\hat{\pi} V \\
@. \widehat\pic(\A^{\Gamma}) @>>> \widehat\pic(G) \\
\end{CD}
\end{equation*}
et on définit l'homomorphisme
$\hat{\psi}:\A^{t,\Gamma'}(S)\rightarrow\widehat\pic(G)$
(associé à
l'inclusion $G\subseteq\A^{\Gamma}$)
comme étant le composé de ces morphismes. Il est clair que $\hat{\psi}$
est un
relèvement de
$\psi$.


\section{Applications}

En guise d'applications de notre construction, nous allons à présent
démontrer les
théorèmes
\ref{intro1} et \ref{intro2} de
l'introduction.

\subsection{Résultat d'annulation sur les points de torsion}
\label{annulsurtors}
Le résultat suivant est une reformulation du théorème 4.1 de \cite{gil}.
Signalons au
passage que
sa démonstration utilise l'hypothèse
d'excellence de $R$.

\begin{thm}
\label{msim}
Supposons que $\E\rightarrow S$ soit le modèle de Néron d'une courbe
elliptique à
réduction
semi-stable sur $K$. Soit $n$ un entier
premier à $6$. Alors la restriction du $\gm$-torseur $t(W)$ au sous-groupe
$\Eo[n]\times_S
\E^t[n]$
est triviale.
\end{thm}

Nous pouvons en déduire le corollaire suivant (généralisant les théorèmes
d'annulation
précédemment
obtenus) :

\begin{cor}
\label{cor}
Soit $\E$ soit le modèle de Néron d'une courbe elliptique à réduction
semi-stable sur
$K$. Soit
$G$
(resp. $H$) un sous-groupe
quasi-fini et fermé de $\Eo$ (resp. $\E$). Soit
$$
\begin{CD}
@.\psi :\E^t(S) @>>> \pic(G) \\
\text{resp. } @.\psi':\E^{t,\circ}(S) @>>> \pic(H) \\
\end{CD}
$$
le morphisme associé à l'inclusion $G\subseteq\Eo$ (resp. $H\subseteq\E$).
Si l'ordre de
$G_{\eta}$
(resp. de $H_{\eta}$) est premier à
$6$, alors $\psi$ (resp. $\psi'$) s'annule sur les points de torsion.
\end{cor}

\begin{proof}
Soit $p\in\E^t(S)$ un point de $m$-torsion. On peut écrire $m=nN$ où $N$
est premier
à l'ordre de $G_{\eta}$, et l'ensemble des facteurs premiers de $n$ est un
sous-ensemble
de
l'ensemble des
facteurs premiers de l'ordre de $G_{\eta}$. Alors $Np$ est un point de
$n$-torsion et,
d'après la
proposition
\ref{oodeux}, il suffit de montrer la
nullité de $\dai(Np)|_G$ pour en déduire celle de $\dai(p)|_G$.

D'autre part, $G$ est un sous-groupe de $\Eo[n]$ (car $G$ est tué par
l'ordre de
$G_{\eta}$), et
$Np$ se factorise à travers $\E^t[n]$.
On peut alors écrire (cf. la remarque \ref{ddef})
$$\dai(Np)|_G=(i\times Np)^*(t(W)|_{\Eo[n]\times_S \E^t[n]})\,.$$
Ainsi, pour montrer que $\dai(p)|_G$ est nul, il suffit de montrer que la
restriction de
$t(W)$ à
$\Eo[n]\times_S \E^t[n]$ est triviale,
ce qui est bien le cas d'après le théorème \ref{msim}. On en déduit le
résultat
d'annulation de
$\psi$, sachant que $\psi(p)=\dai(p)|_G$
pour tout $p\in\E^t(S)$.
Le résultat d'annulation de $\psi'$ se démontre de la même façon.
\end{proof}

Il est clair que le corollaire \ref{cor} implique le théorème
\ref{intro1}.

\begin{exm}
\label{subtil}
Soit $E_L$ une courbe elliptique semi-stable sur un corps de nombres $L$,
on suppose que
$E_L$ a
réduction torique en au moins une place
$v_L$ de $\mathcal{O}_L$. On fixe un nombre premier $N\neq 2,3$ ne
divisant pas l'ordre
du groupe
des composantes de la fibre en $v_L$
du modèle de Néron de $E_L$, et tel que $v_L$ ne divise pas $N$.

Soit $K=L(E_L[N])$, soit $\E\rightarrow S=Spec(\mathcal{O}_K)$ le modèle
de Néron de la
$K$-courbe
elliptique $E_L\times_L K$, et soit
$v$ une place de $K$ divisant $v_L$. Alors $\E[N]$ est fini et plat de
rang $N^2$ sur
$S$, donc
le
groupe $\Phi_v$ des composantes de
$\E_v$ est cyclique d'ordre divisible par $N$, d'après la remarque
\ref{boxall}.

Soit $j$ le $j$-invariant de la courbe $E_L$. Alors l'ordre du groupe des
composantes de
la fibre
en $v_L$ du modèle de Néron de $E_L$
est égal à $-v_L(j)$ (ici, la valuation $v_L$ est normalisée de façon à
prendre la
valeur $1$ sur
une uniformisante).
De même, l'ordre du groupe des composantes de la fibre de $\mathcal{E}$ en
$v$ est égal
à
$-v(j)$.
De plus, $v(j)=e(K/L)v_L(j)$
où $e(K/L)$ est l'indice de ramification de $K/L$. Sachant que $e(K/L)$
divise l'ordre
du groupe
de
Galois de $K/L$, lequel est
isomorphe à un sous-groupe
de $\mathbf{GL}_2(\mathbb{Z}/N\mathbb{Z})$, on en déduit que $N^2$ ne
divise pas
$e(K/L)$.

Ainsi le groupe $\Phi_v$ est cyclique d'ordre $Nk$, avec $k$ premier à
$N$. On peut
alors choisir
un point $p_1\in \E(S)$ d'ordre $N$
qui ne se réduise pas dans $\Eo_v(k_v)$. Il est alors clair que le point
$p_1$
n'appartient pas à
$\E^{N\Phi}(S)$.

Rappelons que, par auto-dualité des courbes elliptiques, on dispose d'un
isomorphisme
canonique
$\E\simeq \E^t$. Soit
$\psi_N:\E(S)\rightarrow \pic(\Eo[N])$ l'homomorphisme associé à
l'inclusion
$\Eo[N]\subseteq
\Eo$.
Considérons la suite exacte
\[
\begin{CD}
0 @>>> \E[N] @>>> \E @>[N]>> \E^{N\Phi} @>>> 0\,.\\
\end{CD}
\]
D'après le paragraphe \ref{kumm}, la restriction de l'homomorphisme
$\psi_N$ à
$\E^{N\Phi}(S)$
étudie la structure galoisienne des
$\E[N]$-torseurs obtenus grâce au cobord de cette suite exacte.

Cependant, on constate ici que notre construction de $\psi_N$ permet
également l'étude
de
torseurs
qui ne proviennent pas du cobord
cette suite exacte. En effet, le point $p_1$ donne naissance à un
$\Hom_S(\Eo[N],\gm)$-torseur
$\delta(p_1)$, dont la structure
galoisienne est triviale d'après le corollaire \ref{cor}, et pourtant
$p_1\not\in\E^{N\Phi}(S)$.
\end{exm}

\subsection{Résultat d'injectivité}
\label{inj}
On suppose à présent que $K$ est un corps de nombres, et que $R$ est
l'anneau des
entiers de $K$.
 
Dans la suite, le groupe $\widehat{\pic}$ sera relatif à l'ensemble
$\mathfrak{S}$ formé
des
places
de mauvaise réduction de $\A$ ainsi
que les places infinies de $K$.

Soit $\ell$ un nombre premier. Alors, pour tout entier $n$, nous disposons
d'un
homomorphisme
$\hat{\psi}_{\ell^n}:\A^t(S)\rightarrow\widehat{\pic}(\Ao[\ell^n])$.
D'autre part,
l'inclusion
$\Ao[\ell^n]\rightarrow \Ao[\ell^{n+1}]$
donne lieu à une flèche de restriction
$\widehat{\pic}(\Ao[\ell^{n+1}])\rightarrow
\widehat{\pic}(\Ao[\ell^n])$, qui fait commuter le
diagramme
\[
\begin{CD}
\A^t(S) @>\hat{\psi}_{\ell^{n+1}}>> \widehat{\pic}(\Ao[\ell^{n+1}]) \\
@VVV @VVV \\
\A^t(S) @>\hat{\psi}_{\ell^n}>> \widehat{\pic}(\Ao[\ell^n])\,.\\
\end{CD}
\]
On obtient ainsi, par passage à la limite (projective), un homomorphisme
$$\hat{\Psi}_{\ell}=\limp\hat{\psi}_{\ell^n}:\A^t(S)\otimes\mathbf{Z}_{\ell}\longrightarrow
\limp
\widehat{\pic}(\Ao[\ell^n])\,.$$

On peut de même définir, pour le $U$-schéma abélien $A_U$, un morphisme
analogue, que
nous
noterons
$\hat{\Psi}_{U,\ell}$. Ce dernier
coïncide avec le morphisme défini par Agboola et Pappas dans \cite{ap}.
Nous pouvons
énoncer pour
$\hat{\Psi}_{U,\ell}$ le résultat
suivant (voir \cite{ap}, Theorem 6.4 et Theorem 1.2) :

\begin{thm}(Agboola-Pappas).\label{agpa}
L'homomorphisme $\hat{\Psi}_{U,\ell}$ est injectif modulo les points de
torsion. En
outre, si
tous les points de $S$ de caractéristique $\ell$ sont contenus dans $U$,
et si $\ell$ ne
divise
pas $6\cdot\disc(K/\mathbb{Q})$, alors $\hat{\Psi}_{U,\ell}$ est injectif.
\end{thm}

\begin{proof}[Démonstration du théorème \ref{intro2}]
Nous avons un diagramme commutatif
\[
\begin{CD}
\A^t(S) @>\hat{\psi}_{\ell^n}>> \widehat{\pic}(\Ao[\ell^n]) \\
@VVV @Vj^*VV \\
A_U^t(U) @>\hat{\psi}_{U,\ell^n}>> \widehat{\pic}(A_U[\ell^n])\,.\\
\end{CD}
\]
Ainsi, par passage à la limite, la composée de $\hat{\Psi}_{\ell}$ avec le
morphisme
naturel
$$\limp \widehat{\pic}(\Ao[\ell^n])\longrightarrow \limp
\widehat{\pic}(A_U[\ell^n])$$
est égale au morphisme $\hat{\Psi}_{U,\ell}$ défini plus haut. En
particulier, tout
résultat
d'injectivité portant sur
$\hat{\Psi}_{U,\ell}$ donne lieu au même résultat pour
$\hat{\Psi}_{\ell}$. Ainsi le
théorème
\ref{agpa} entraîne le théorème
\ref{intro2}.
\end{proof}



\begin{thebibliography}{12}

\bibitem[A1]{a1} \textsc{A. Agboola}, {\it A geometric description of the
class
invariant
homomorphism}, J. Théor. Nombres Bordeaux {\bf
6} (1994), 273--280.
\bibitem[A2]{a2} ---------, {\it Torsion points on elliptic curves and
Galois module
structure},
Invent. Math. {\bf 123} (1996),
105--122.
\bibitem[A-P]{ap} \textsc{A. Agboola} and \textsc{G. Pappas}, {\it On
arithmetic class
invariants},
Math. Ann. {\bf 320} (2001),
339--365.
\bibitem[A-T]{at} \textsc{A. Agboola} and \textsc{M. J. Taylor}, {\it
Class invariants of Mordell-Weil groups}, J. Reine Angew. Math. {\bf 447}
(1994), 23-61.
\bibitem[An]{anan} \textsc{S. Anantharaman}, {\it Schémas en groupes,
espaces homogènes
et
espaces
algébriques sur une base de dimension
$1$}, Bull. Soc. Math. Fr., Suppl., Mém. {\bf 33} (1973), 5--79.
\bibitem[BLR]{ray} \textsc{S. Bosch}, \textsc{W. Lütkebohmert} and
\textsc{M. Raynaud},
{\it
Néron
Models}, Ergeb. Math. Grenzgeb. (3), vol. 21 (Springer, Berlin-Heidelberg-New York, 1990).
\bibitem[CN-T]{cnt} \textsc{P. Cassou-Noguès} et \textsc{M. J. Taylor},
{\it Structures
galoisiennes et courbes elliptiques}, J. Théor.
Nombres Bordeaux {\bf 7} (1995), 307--331.
\bibitem[G]{gil} \textsc{J. Gillibert}, {\it Invariants de classes : le
cas
semi-stable}, Compositio Mathematica {\bf 141} (2005), 887--901.
\bibitem[SGA 4]{gro4} \textsc{A. Grothendieck}, \textsc{M. Artin} et
\textsc{J. L.
Verdier}, {\it
Théorie des topos et cohomologie étale
des schémas}, Lecture Notes in Mathematics, vols 269, 270 (Springer, Berlin-Heidelberg-New York, 1972).
\bibitem[SGA 7]{gro7} \textsc{A. Grothendieck}, {\it Groupes de monodromie
en géométrie
algébrique}, Lecture Notes in Mathematics, vol. 288 (Springer, Berlin-Heidelberg-New York, 1972).
\bibitem[Mi]{mil} \textsc{J. S. Milne}, {\it Arithmetic Duality Theorems},
Perspectives
in
Mathematics, vol. 1 (Academic Press, Boston, MA, 1986).
\bibitem[Mu]{mu} \textsc{D. Mumford}, {\it Bi-extensions of formal
groups}, in {\it Algebraic Geometry} (Bombay, 1968), (Oxford University Press, 1969), 307--322.
\bibitem[P]{p1} \textsc{G. Pappas}, {\it On torsion line bundles and
torsion points on
abelian
varieties}, Duke Math. J. {\bf 91}
(1998), 215--224.
\bibitem[S-T]{st} \textsc{A. Srivastav} and \textsc{M. J. Taylor}, {\it
Elliptic curves
with
complex multiplication and Galois module
structure}, Invent. Math. {\bf 99} (1990), 165--184.
\bibitem[T88]{t2} \textsc{M. J. Taylor}, {\it Mordell-Weil groups and the
Galois module
structure
of rings of integers}, Illinois J.
Math. {\bf 32} (1988), 428--452.
\bibitem[T91]{t3} \textsc{M. J. Taylor}, {\it $L$-functions and Galois
modules :
Explicit Galois
Modules}, in {\it $L$-functions and
Arithmetic}, LMS Lecture Notes, vol. 153 (Cambridge University Press,
1991).
\bibitem[W]{w} \textsc{W. C. Waterhouse}, {\it Principal homogeneous
spaces and group
scheme
extensions}, Trans. Am. Math. Soc. {\bf
153} (1971),
181--189.
\bibitem[We]{werner} \textsc{A. Werner}, {\it On Grothendieck's pairing of
component
groups in
the
semistable reduction case}, J. Reine
Angew. Math. {\bf 486} (1997), 205--215.

\end{thebibliography}
\end{document}